\documentclass[preprint,12pt]{elsarticle}

\usepackage{amssymb}

\journal{Linear Algebra and its Applications, 435(2011) 2784-2792}

\begin{document}

\begin{frontmatter}



\title{An algorithm for determining copositive matrices\tnoteref{t1}}
\tnotetext[t1]{This research was supported by the National Natural
Science Foundation of China(11001228,10901116). }


\author[jx]{Jia Xu\corref{cor1}}
\ead{j.jia.xu@gmail.com}

\author[yy]{Yong Yao}
\ead{yaoyong@casit.ac.cn}

\cortext[cor1]{Corresponding author.}

\address[jx]{College of Computer Science and Technology,
Southwest University for Nationalities, Chengdu, Sichuan 610041,
China}

\address[yy]{Chengdu Institute of Computer Applications, Chinese
Academy of Sciences, Chengdu, Sichuan 610041, China}

\begin{abstract}
In this paper, we present an algorithm of simple exponential growth
called COPOMATRIX for determining the copositivity of a real
symmetric matrix. The core of this algorithm is a decomposition
theorem, which is used to deal with simplicial subdivision of
$\widehat{T}^{-}=\{y\in \Delta_{m}| \beta^Ty\leq 0\}$ on the
standard simplex $\Delta_m$, where each component of the vector
$\beta$ is -1, 0 or 1.

\end{abstract}


\end{frontmatter}


\section{Introduction}

\textbf{QUESTION 1}\quad  Let $A$ be a given $n\times n$ real
symmetric matrix, ${\Bbb R}_+$ be the set of nonnegative real
numbers, and
$$
Q(X)=X^TAX,\  X\neq 0
$$ be a quadratic form. What conditions should $A$ satisfy for
$[\forall X\in {\Bbb R}^n_+ ,\ Q(X)\geq 0\ (>0)]$?

If $[\forall X\in {\Bbb R}^n_+ ,\ Q(X)\geq 0\ (>0)]$, then the
quadratic form $Q(X)$ is called a (strictly) copositive quadratic
form and the corresponding matrix $A$ is called a (strictly)
copositive matrix.

Copositive matrices have numerous applications in diverse fields of
applied mathematics, especially in mathematical programming and
graph theory (see
\cite{Bomze:00,Bomze:02,Burer,Danninger:00,Eichfelder,Hadjicostas,Hall,
Juttler,Klerk,Lemke,Quist}). Therefore copositivity has been studied
thoroughly since 1950s (see
\cite{Andersson,Bundfuss,Hadeler,Haynsworth,Horn,Ikramov:01,Johnson:00,
Kaplan:00,Kaplan:01,Martin:00,Martin:01,Motzkin,Valiaho:00,
Valiaho:01,Valiaho:02,S.Yang:00,S.Yang:01}).

In general, it is an NP-complete problem to determine whether a
given $n\times n$ symmetric matrix is not copositive
\cite{Murty,Parrilo}. This means that every algorithm that solves
the problem, in the worst case, will require at least an exponential
number of operations, unless P=NP. For that reason, it is still
valuable for the existence of so many incomplete algorithms
discussing some special kinds of matrices (see
\cite{Bomze:00,Bomze:01,Cottle,Danninger:01,Hoffman,Hogben,Ikramov:00,Johnson:01,Li,Parrilo}).
For small values of $n$($\leq 6$), some necessary and sufficient
conditions have been constructed (see
\cite{Andersson,Hadeler,Li,S.Yang:01}). From another viewpoint,
\textbf{QUESTION 1} is a typical real quantifier elimination problem
\cite{Basu,Collins:00,Collins:01,Marshall,Mishra,Tarski,Yang:00,Yang:01},
which can be solved by standard tools of real quantifier elimination
(e.g., using CAD) \cite{Basu,Collins:00,Collins:01,Yang:02,Yang:03}.
Thus, there is a complete algorithm for determining copositive
matrices theoretically. Unfortunately, this algorithm is not
efficient in practice for the CAD algorithm is of doubly exponential
time complexity (see \cite{Basu,Collins:00,Collins:01}). In this
paper, we will construct a complete algorithm with singly
exponential time bound.

The standard simplex $\Delta_m  (m\geq 2)$ is defined as the
following set
$$\Delta_m=\{(y_1,\ldots,y_m)^T|\ y_1+\cdots +y_m=1, y_1\geq 0,\ldots,y_m\geq 0\}.$$
It is well known that the dimension of $\Delta_m$ is $m-1$. Denote
the vertices of $\Delta_m$ as $e_1,\ldots,e_m$, namely,
$e_1=(1,0,\ldots,0)^T,\ldots,e_m=(0,0,\ldots,1)^T$.

Let $A\in {\Bbb R}^{n\times n}$ be symmetric and be partitioned as
$$
A=[\alpha_{ij}]=\left[
\begin{array}{cc}
\alpha_{11}& \alpha^T \\
\alpha& A_2
\end{array}
\right].
$$
Define $B=\alpha_{11}A_2-\alpha \alpha^T$. It is easy to see the
following facts (cf.\cite{Andersson})

1. If $\alpha_{1i}\geq 0,\ i=2,\ldots,n$, then $A$ is (strictly)
copositive $\Longleftrightarrow \alpha_{11}\geq 0 \ (>0)$ and $A_2$
is (strictly) copositive.

2. If at least one of $\alpha_{1i}$ is negative, then we need only
to focus on the set of points $T^{-}=\{y\in \Delta_{n-1} |\
\alpha^Ty \leq 0\}$. It is well known that $T^{-}$ is a convex
polytope on $\Delta_{n-1}$ (see \cite{Andersson}). The polytope
$T^{-}$ can be subdivided into the simplices $S_1,\ldots,S_p$, that
is,
$$T^{-}=\bigcup_{i=1}^p S_i,\ {\rm int} (S_i)\bigcap {\rm int} (S_j)=
\emptyset ,{\rm for}\ i\neq j,$$ where ${\rm int} (S_i)$ denotes the
interior of simplex $S_i$. The coordinates of the vertices that span
the simplex $S_i$ constitute a matrix denoted as $W_i$. Andersson et
al. (\cite{Andersson},\ p.23) proved the following results.

\newtheorem{lem1}{Lemma}[section]
\begin{lem1}\upshape\label{lem1}
(a) $A$ is copositive iff $\alpha_{11}\geq 0$ and $A_2,\
W_{1}^TBW_1,\ \ldots,\  W_{p}^TBW_p$ are all copositive.

(b) $A$ is strictly copositive iff $\alpha_{11}> 0$ and $A_2,\
W_{1}^TBW_1,\ \ldots,\  W_{p}^TBW_p$ are all strictly copositive.
\end{lem1}

In order to formulate the algorithm of Lemma \ref{lem1}, we first
consider how to obtain the simplicial subdivision of the polytope
$T^{-}=\{y\in \Delta_{n-1} |\ \alpha^Ty\leq 0\}$. For small values
of $n(\leq 6)$,\ Andersson et al.\cite{Andersson} and Yang and Li
\cite{S.Yang:01} give the simplicial subdivision of $T^{-}$.
However, they do not provide a procedure for a simplicial
subdivision of $T^{-}$ for arbitrary values of  n.
 We propose a simplicial subdivision of $T^{-}$ for
all values of $n$, and consequently construct a complete algorithm
for determining the copositivity of an $n\times n$ matrix.

We will adopt a flexible approach. Rather than subdivide $T^{-}$
into simplices (of course our method is also valid for subdividing
$T^{-}$ into simplices), we first transform the matrix $A$ into the
following matrix called $\widehat{A}$.

Let $\alpha=(\alpha_{12},\ldots,\alpha_{1n})^T$ and $D={\rm
diag}(d_{2},\ldots,d_{n}),$ where
$$
d_{i}=\left\{
\begin{array}{ll}
1,& {\rm if}\ \alpha_{1i}=0;\\
1/|\alpha_{1i}|,& {\rm if}\ a_{1i}\neq 0.
\end{array}
\right .
$$
Then
\begin{equation}\label{eqa0}
\widehat{A}= \left[
\begin{array}{cc}
1& 0 \\
0& D
\end{array}
\right]A \left[
\begin{array}{cc}
1& 0 \\
0& D
\end{array}
\right] =\left[
\begin{array}{cc}
\alpha_{11}& \widehat{\alpha}^T \\
\widehat{\alpha}& DA_2D
\end{array}
\right].
\end{equation}
where $\widehat{\alpha}=({\rm sign}(\alpha_{12}),\ldots,{\rm
sign}(\alpha_{1n}))^T$.

Obviously, $A$ is (strictly) copositive
$\Longleftrightarrow$$\widehat{A}$ is (strictly) copositive. Apply
Lemma \ref{lem1} to $\widehat{A}$. Let
$$\beta_{1}={\rm sign}(\alpha_{12}),\ldots,\ \beta_{n-1}={\rm sign}(\alpha_{1n}).$$
Thus we just need to subdivide $\widehat{T}^-$ into simplices, where
$$\widehat{T}^-=\{y\in \Delta_{n-1} |(\beta_{1},\ldots,\beta_{n-1})y\leq
0,\ \beta_{i}\in \{-1,0,1\}\}.$$

Next we make further simplification: separate -1,0,1 from
$\beta_{1},\ldots,\beta_{n-1}$, namely let
$$
\begin{array}{c}
\beta_{a_1}=\cdots=\beta_{a_s}=1,\
\beta_{b_1}=\cdots=\beta_{b_t}=-1,\
\beta_{c_1}=\cdots=\beta_{c_r}=0.\\
\{a_1,\ldots,a_s,\ b_1,\ldots,b_t,\ c_1,\ldots,\ c_r \}=
\{1,\ldots,n-1\},\\
 r,s,t\geq 0,\ t\geq 1,\ \ r+s+t=n-1.
\end{array}
 $$
In geometry it is easy to see that the convex polytope
$\widehat{T}^-$ is the convex hull of its surface $S^-$ and its
vertices $e_{c_1},\ldots,e_{c_r}$, that is,
\begin{equation}\label{eqa1}
\widehat{T}^-={\rm conv}\{e_{c_1},\ldots,e_{c_r},S^-\}.
\end{equation}
\begin{equation}\label{eqa2}
\begin{array}{l}
S^-=\{(y_1,\ldots,y_{n-1})^T\in \Delta_{n-1}|y_{a_1}+\cdots+y_{a_s}-y_{b_1}-\cdots-y_{b_t}\leq 0,\\
\qquad \quad (y_{a_1},\ldots,y_{a_s},y_{b_1},\ldots,y_{b_t})^T\in
\Delta_{s+t}\}.
\end{array}
\end{equation}

If the simplicial subdivision of $S^-$ is known, the simplicial
subdivision of $\widehat{T}^-$ is directly obtained by (\ref{eqa1}).
So we just need to study the simplicial subdivision of the polytope
$S^-$.

\section{A simplicial subdivision algorithm for the convex polytope $S^-$}

\subsection{Fundamental notations}
The notation $S^-$ is simple, but it can not reveal the information
of convex polytopes. In order to simplify the descriptions, we will
introduce a new notation, which is fundamental to our study.

\newdefinition{def1}{Definition}[section]
\begin{def1}\upshape\label{def1}
Suppose that two sequences of positive integers
$[a_1,a_2,\ldots,a_s]$,\\$ [b_1,b_2,\ldots,b_t]$ satisfy
$$\{a_1,\ldots,a_s,\ b_1,\ldots,b_t\}\subseteq \{1,2,\ldots,m\},\
 s\geq 0,t\geq 1,m\geq s+t\geq 2,$$
where all of $s+t$ elements of $\{a_1,\ldots,a_s,\ b_1,\ldots,b_t\}$
are distinct. Then the notation
$[[a_1,a_2,\ldots,a_s],[b_1,b_2,\ldots,b_t]]_m$ is defined as the
polytope $S^-$ in (\ref{eqa2}).
\end{def1}

For example, let us compare the polytope $[[2,3],[5]]_5$ and the
polytope $[[2,3],[5]]_6$. \ $[[2,3],[5]]_5$ denotes the polytope
$$\{(y_1,\ldots,y_5)^T\in \Delta_5 |y_2+y_3-y_5\leq 0,(y_2,y_3,y_5)^T\in \Delta_3\}.$$
Here $(y_2,y_3,y_5)^T\in \Delta_3$ implies that $y_1=0,y_4=0$. \
$[[2,3],[5]]_6$ indicates the polytope
$$\{(y_1,\ldots,y_6)^T\in \Delta_6 |y_2+y_3-y_5\leq 0,(y_2,y_3,y_5)^T\in \Delta_3\}.$$
Here $(y_2,y_3,y_5)^T\in \Delta_3$ implies that $y_1=0,y_4=0,y_6=0$.
It is clear that $[[2,3],[5]]_5$ and $[[2,3],[5]]_6$ are congruent,
although they are sets in simplices of different dimensions.

For $0\leq k \leq m-1$, the polytope $L_k^-$ is defined as
$$L_k^-=\{(y_1,\ldots,y_m)^T\in \Delta_m|y_1+\cdots+y_k-y_{k+1}-\cdots-
y_{m}\leq 0\}.$$

$L_k^-$ is written as $[[1,\ldots,k],[k+1,\ldots,m]]_m$ by the
notation of Definition \ref{def1}. $L_k^-$ is a special case of
$S^-$, but this notation is more convenient for our analysis.

In the following we will study the basic geometric properties of
convex polytope
$[[a_1,a_2,\ldots,a_s],[b_1,b_2,\ldots,b_t]]_m$.\\

\subsection{Geometric properties of $S^-$}

Let $e_1,\ldots,e_m$ be vertices of the standard simplex $\Delta_m$,
and $M_{i,j}=(e_i+e_j)/2$ be the midpoint of the line segment
$\overline{e_ie_j}$.

The following result is stated in \cite{Andersson} without proof.
For completeness, we give a proof.

\newtheorem{lem2}{Lemma}[section]
\begin{lem2}\upshape\textbf{\cite{Andersson}}\label{lem2}
\ \ Given a convex polytope $L_k^{-}$, then all of its vertices are
$$V=\{e_{k+1},\ldots,e_m,\ M_{i,j},i=1,2,\ldots,k,\ j=k+1,\ldots,m\}.$$
The number of the vertices is $|V|=(k+1)(m-k)$.
\end{lem2}

\newproof{pf2}{Proof}
\begin{pf2}
Note that the convex polytope $L_k^{-}$ is obtained by cutting the
standard simplex $\Delta_m$ with the hyperplane
$$L_{=0}:\ y_1+\cdots+ y_k- y_{k+1}-\cdots-y_{m}=0.$$
Therefore the vertices of the polytope $L_k^-$ come from two parts:
one part is vertices of $\Delta_m$, that is,
$\{e_{k+1},\ldots,e_m\}$; while the other part is the intersection
points of the hyperplane $L_{=0}$ and the edges of standard simplex
$\Delta_m$.

First consider the intersection point of $L_{=0}$ and the edge
$ae_1+be_{k+1}\ (a,b\geq 0,\ a+b=1)$. Substitute $ae_1+be_{k+1}$
into the following equations,
$$y_1+\cdots+ y_k- y_{k+1}-\cdots-y_{m}=0,\ y_1+\cdots+y_m=1$$
Therefore, the solutions are $a=1/2,b=1/2$, namely, the intersection
point is $M_{1,k+1}$.

In the same way, we get all intersection points of $L_{=0}$ and the
edges of $\Delta_m$. They are $\{M_{i,j},i=1,2,\ldots,k,\
j=k+1,\ldots,m.\}$.

Hence the number of vertices of $L_k^-$ is
$|V|=m-k+k(m-k)=(k+1)(m-k)$.
\end{pf2}

Likewise, we can prove the following lemma.

\newtheorem{lem3}[lem2]{Lemma}
\begin{lem3}\upshape\label{lem3}
Given a convex polytope
$[[a_1,a_2,\ldots,a_s],[b_1,b_2,\ldots,b_t]]_m$, then all of its
vertices are
$$V=\{e_{b_1},\ldots,e_{b_t},\ M_{a_i,b_j},i=1,2,\ldots,s,\ j=1,\ldots,t\}.$$
The number of the vertices is $|V|=(s+1)t$.
\end{lem3}

We see that the polytope
$[[a_1,a_2,\ldots,a_s],[b_1,b_2,\ldots,b_t]]_m$ and the polytope
$L_k^{-}$ are similar in many respects, which will be further
discussed.

\newtheorem{lem4}[lem2]{Lemma}
\begin{lem4}\upshape\label{lem4}
The convex polytope $L_k^{-}$ is simplicial iff $k=0$, or $k=m-1$.
\end{lem4}

\newproof{pf4}{Proof}
\begin{pf4}
When $k=0$, $L_0^-=\Delta_m$ is simplicial.

When $k=m-1$, consider the convex polytope
$$
L_{m-1}^{-}:=\{(y_1,\ldots,y_m)^T \in \Delta_m| \ y_1+\cdots+
 y_{m-1}-y_{m}\leq 0\}.
$$
By Lemma \ref{lem2}, we know that all vertices of $L_{m-1}^{-}$ are
$\{e_m,\ M_{i,m},i=1,2,\ldots,$\\$m-1\}$. Obviously all the vectors
of $\{M_{i,m}-e_m,i=1,2,\ldots,m-1 \}$ are linearly independent, so
$L_{m-1}^{-}$ is simplicial.
\end{pf4}

Conversely, we know that the dimension of the polytope $L_k^{-}$ is
$m-1$. If $k\neq 0,m-1$, then by Lemma \ref{lem2}, the number of the
vertices of $L_k^{-}$ is $(k+1)(m-k)\neq m$, so $L_k^{-}$ is not
simplicial.

\newtheorem{lem5}[lem2]{Lemma}
\begin{lem5}\upshape\label{lem5}
The convex polytope $[[a_1,a_2,\ldots,a_s],[b_1,b_2,\ldots,b_t]]_m$
(here the vertices are obtained by Lemma \ref{lem3}) is simplicial
iff $s=0$, or $t=1$.
\end{lem5}

\newtheorem{lem6}[lem2]{Lemma}
\begin{lem6}\upshape\label{lem6}
The dimension of the polytope
$[[a_1,a_2,\ldots,a_s],[b_1,b_2,\ldots,b_t]]_m$ is $(s+t-1)$.
\end{lem6}

If the polytope $[[a_1,a_2,\ldots,a_s],[b_1,b_2,\ldots,b_t]]_m$ is
not a simplex, we will subdivide it into simplices.

\newtheorem{lem7}[lem2]{Lemma}
\begin{lem7}\upshape\label{lem7}
If the polytope $[[a_1,a_2,\ldots,a_s],[b_1,b_2,\ldots,b_t]]_m$ is
not a simplex, then there are only two $(s+t-2)$-dimensional
surfaces that do not include the vertex $M_{a_1,b_1}$. They are
$$
[[a_2,\ldots,a_s],[b_1,b_2,\ldots,b_t]]_m,\
[[a_1,\ldots,a_s],[b_2,\ldots,b_t]]_m.\
$$
(obtained by deleting $a_1,b_1$ from array
$[[a_1,a_2,\ldots,a_s],[b_1,b_2,\ldots,b_t]]_m$ respectively)
\end{lem7}

\newproof{pf7}{Proof}
\begin{pf7}
All the $(s+t-2)$-dimensional surfaces of the convex polytope\\
$[[a_1,a_2,\ldots,a_s],\ [b_1,b_2,\ldots,b_t]]_m$ are obviously

\quad \quad \quad \quad $[[\widehat{a}_1,
a_2,\ldots,a_s],[b_1,b_2,\ldots,b_t]]_m,$

\quad \quad \quad \quad
$[[a_1,\widehat{a}_2,\ldots,a_s],[b_1,b_2,\ldots,b_t]]_m,$

\quad \quad \quad \quad $\ldots,$

\quad \quad \quad \quad
$[[a_1,a_3,\ldots,a_s],[b_1,b_2,\ldots,\widehat{b}_t]]_m$
\\
(where the notation $[[\widehat{a}_1,
a_2,\ldots,a_s],[b_1,b_2,\ldots,b_t]]_m$ is the polytope with $a_1$
removed) and
$$
\begin{array}{l}
\{(y_1,\ldots,y_m)^T\in
\Delta_m|y_{a_1}+\cdots+y_{a_s}-y_{b_1}-\cdots-y_{b_t}= 0,
\ (y_{a_1},\ldots,y_{a_s},\\
\quad y_{b_1},\ldots,y_{b_t})^T\in \Delta_{s+t}\}.
\end{array}
$$
That makes $s+t+1$ $(s+t-2)$-dimensional surfaces in all. We can
verify that only
$$
[[a_2,\ldots,a_s],[b_1,b_2,\ldots,b_t]]_m,\
[[a_1,\ldots,a_s],[b_2,\ldots,b_t]]_m
$$
do not include the vertex $M_{a_1,b_1}$.
\end{pf7}

Lemma \ref{lem7} leads to the following decomposition theorem.\\

\subsection{The decomposition process for the polytope $S^-$}

\newtheorem{thm}{Theorem}[section]
\begin{thm}[decomposition theorem]\upshape\label{thm}
If the polytope $[[a_1,a_2,\ldots,a_s]$,\\ $[b_1,b_2,\ldots,b_t]]_m$
is not simplicial, then it can be decomposed into the union of two
convex polytopes (not always simplicial). The expression is
$$
\begin{array}{l}
\quad [[a_1,a_2,\ldots,a_s],[b_1,b_2,\ldots,b_t]]_m\\
 ={\rm
 conv}\{M_{a_1,b_1},[[a_2,\ldots,a_s],[b_1,b_2,\ldots,b_t]]_m\}\bigcup \\
 \quad {\rm conv}\{M_{a_1,b_1},[[a_1,a_2,\ldots,a_s],
 [b_2,\ldots,b_t]]_m\}.
\end{array}
$$
Here ${\rm conv}\{S\}$ denotes the convex hull of the set $S$ of
points .
\end{thm}

\newproof{pf}{Proof}
\begin{pf}
This follows from Lemma \ref{lem7}.
\end{pf}

Based on Theorem \ref{thm}, the polytope $S^-$ can be easily
subdivided into simplices.

\newtheorem{eg1}{Example}
\begin{eg1}\upshape\label{eg1}
Show the simplicial subdivision of the following convex polytope
$$L_2^- :=\{(y_1,\ldots,y_5)|\ y_1+y_2-y_3-y_4-y_5\leq 0,(y_1,\ldots,y_5)^T\in \Delta_5\}.$$
\end{eg1}

\newproof{st1}{Solution}
\begin{st1}
Denote $L_2^-$ as $[[1,2],[3,4,5]]_5$. We know that
$[[1,2],[3,4,5]]_5$ is not simplicial by Lemma \ref{lem5}. Using
Theorem \ref{thm} we have
$$
\begin{array}{l}
[[1,2],[3,4,5]]_5 ={\rm conv}\{M_{1,3},[[2],[3,4,5]]_5\}\bigcup {\rm
conv}\{M_{1,3},[[1,2],[4,5]]_5\}.
\end{array}
$$

By Lemma \ref{lem5} we know that both $[[2],[3,4,5]]_5$ and
$[[1,2],[4,5]]_5$ are not simplicial. Therefore we repeatedly apply
Theorem \ref{thm} to them  and have\\

$
\begin{array}{l}
[[2],[3,4,5]]_5\\
={\rm conv}\{M_{2,3},[[\ ],[3,4,5]]_5\}\bigcup
{\rm conv}\{M_{2,3},[[2],[4,5]]_5\}.\\
={\rm conv}\{M_{2,3},[[\ ],[3,4,5]]_5\}\bigcup {\rm
conv}\{M_{2,3},M_{2,4},[[\ ],[4,5]]_5\}\bigcup\\
\quad {\rm conv}\{M_{2,3},M_{2,4},[[2],[5]]_5\}\\
={\rm conv}\{M_{2,3},e_3,e_4,e_5\}\bigcup {\rm
conv}\{M_{2,3},M_{2,4},e_4,e_5\}\bigcup\\
\quad {\rm conv}\{M_{2,3},M_{2,4},M_{2,5},e_5\}.
\end{array}
$
\\
\\

$
\begin{array}{l}
[[1,2],[4,5]]_5\\
={\rm conv}\{M_{1,4},[[2],[4,5]]_5\}\bigcup
{\rm conv}\{M_{1,4},[[1,2],[5]]_5\}.\\
={\rm conv}\{M_{1,4},M_{2,4},[[\ ],[4,5]]_5\}\bigcup {\rm
conv}\{M_{1,4},M_{2,4},[[2 ],[5]]_5\}\bigcup\\
\quad {\rm conv}\{M_{1,4},[[1,2],[5]]_5\}\\
={\rm conv}\{M_{1,4},M_{2,4},e_4,e_5\}\bigcup {\rm
conv}\{M_{1,4},M_{2,4},M_{2,5},e_5\}\bigcup\\
\quad {\rm conv}\{M_{1,4},M_{1,5},M_{2,5},e_5\}.
\end{array}
$
\\

Finally we get the expression of simplicial subdivision of
$[[1,2],[3,4,5]]_5$,\\

$
\begin{array}{l}
[[1,2],[3,4,5]]_5\\
={\rm conv}\{M_{1,3}, M_{2,3}, e_3, e_4, e_5\}\bigcup
{\rm conv}\{M_{1,3}, M_{2,3}, M_{2,4}, e_4, e_5\}\bigcup\\
\quad {\rm conv}\{M_{1,3},M_{2,3},M_{2,4},M_{2,5},e_5\}
\bigcup {\rm conv}\{M_{1,3},M_{1,4},M_{2,4},e_4,e_5\}\bigcup \\
\quad {\rm conv}\{M_{1,3},M_{1,4},M_{2,4},M_{2,5},e_5\}\bigcup {\rm
conv}\{M_{1,3},M_{1,4},M_{1,5},M_{2,5},e_5\}.
\end{array}
$
\\
\\

So $[[1,2],[3,4,5]]_5$ is a union of six 4-dimensional simplices.
\end{st1}

We summarize the decomposition process of Example \ref{eg1} into the following algorithm.\\

\hrule depth0pt height0.25truemm width\textwidth

\indent{\bf Algorithm 1 (Vmatrix)}

\hrule depth0pt height0.25truemm width\textwidth\quad

Input: The expression of polytope
$[[a_1,a_2,\ldots,a_s],[b_1,b_2,\ldots,b_t]]_m$.

Output: Simplices $D_1,D_2,\ldots,D_p$(denoted by matrices) such
that
$$[[a_1,a_2,\ldots,a_s],[b_1,b_2,\ldots,b_t]]_m=\bigcup_{i=1}^p
D_i,\ {\rm int} (D_i)\bigcap {\rm int} (D_j)=\emptyset ,\ {\rm for}\
i\neq j. $$

V1: Let $F:=\{[[a_1,a_2,\ldots,a_s],[b_1,b_2,\ldots,b_t]]_m\}$,\
${\rm temp}:=\emptyset $.

V2: When $F\neq \emptyset $, repeat the following procedures

\qquad V21: Choose a polytope $N\in F$. If $N$ is simplicial, then

\qquad \qquad \ \ ${\rm temp}:={\rm temp}\cup \{N \}$,
$F:=F\setminus \{N\}$.

\qquad V22: If the polytope $N$ is not simplicial, then by Theorem
\ref{thm} de-

\qquad \qquad \ \ compose it into two convex polytope $B_1,B_2$.

\qquad \qquad \ \  $F:=F\setminus \{N\} \cup \{B_1,B_2\}$. Go to
step V2.

V3: Return ${\rm temp}$.
\\
\hrule depth0pt height0.25truemm width\textwidth \quad

We have written a function in Maple \cite{Maple} to implement the
above algorithm.

Lastly, we will present a formula for computing the number of
simplices given by the polytope $L_k^{-}$ subdivision.

\newtheorem{lem8}[lem2]{Lemma}
\begin{lem8}\upshape\label{lem8}
According to algorithm Vmatrix, the convex polytope
$[[1,\ldots,k]$,\\$[k+1,\ldots,m]]_m$ $(0 \leq k \leq m-1,m\geq 2)$
can be subdivided just into $f(k,m)$ simplices, where
$$f(k,m)=\left(
\begin{array}{c}
m-1\\
k
\end{array}
\right) =\frac{(m-1)!}{k!(m-1-k)!}.
$$
\end{lem8}

We know that $f(k,m)$ has the same recurrence formula as binomial
coefficients by Theorem \ref{thm}. Thus the proof of Lemma
\ref{lem8} is easy via an induction argument.
This formula will be used to estimate the cost of Algorithm 2 in the next section.\\

\section{Determining algorithm for copositive matrices}

In this section, we will present the complete determining algrorithm
of a copositive matrix.

Given an $n\times n$ symmetric matrix
$$
A=[\alpha_{ij}]=\left[
\begin{array}{cc}
\alpha_{11}& \alpha^T \\
\alpha& A_2
\end{array}
\right],
$$

compute $\widehat{A}$\ (see (\ref{eqa0}))
$$
\widehat{A} =\left[
\begin{array}{cc}
\alpha_{11}& \widehat{\alpha}^T \\
\widehat{\alpha}& DA_2D
\end{array}
\right].
$$

Let $B=\alpha_{11}DA_2D-\widehat{\alpha} \widehat{\alpha}^T$, and
let
$$\widehat{\alpha}=({\rm sign}(\alpha_{12}),\ldots,{\rm sign}
(\alpha_{1n}))^T=(\beta_1,\ldots,\beta_{n-1})^T.$$

Define the projection operator \textbf{Proj} of the matrix $A$ as
follows,

$\bullet$ If $\beta_i\geq 0,\ i=1,\ldots,n-1$, then
$${\rm Proj(A)}=\{DA_2D\}.$$

$\bullet$ If there is at least one -1 in $\beta_i$, then
$${\rm Proj(A)}=\{DA_2D,\  W_1^TBW_1,\ \ldots,\  W^T_pBW_P\}.$$
Here the matrices $W_1,\ldots,W_p$ is fixed by the simplicial
subdivision of the convex polytope $\widehat{T}^-$ (see
(\ref{eqa1})).\\

\hrule depth0pt height0.25truemm width\textwidth

\indent{\bf Algorithm 2 (COPOMATRIX)} \hrule depth0pt
height0.25truemm width\textwidth \quad

Input: Symmetric matrix $A\in {\mathbf R}^{n\times n}(n\geq 2)$.

Output: $A$ is copositive, or $A$ is not copositive.

C1: Let $F:=\{A \}$.

C2: Repeat the following steps for the set $F$.

\qquad C21: If the set $F$ is empty, then return ``$A$ is
copositive".

\qquad C22: Check the (1,1)$^{th}$ entry  of every matrix $K$ in set
$F$. If at least

\qquad \qquad \ \ one of them is negative, then return ``$A$ is not
copositive".

\qquad C23: Compute the projective set $P:=\bigcup_{K\in F} {\rm
Proj}(K)$ of set $F$.

\qquad \qquad \ \ $F:=P\setminus \{\hbox{the nonnegative matrices of
$P$} \}$. Go to step C21.
\\
\hrule depth0pt height0.25truemm width\textwidth \quad

Note that the above algorithm is also valid for $2\times 2$
matrices. Furthermore, for strictly copositive matrices we can also
formulate a similar algorithm.

The correctness of the algorithm COPOMATRIX is guaranteed by Lemma
\ref{lem1}, and the algorithm obviously terminates. The cost of the
algorithm mainly depends on the number of simplicial subdivisions of
the polytope. According to Lemma \ref{lem8}, we can estimate that in
the worst case it is at most:
$$
\begin{array}{l}
(\left(
\begin{array}{c}
n-2\\
\left [\frac{n-2}{2}\right ]
\end{array}
\right)+1) (\left(
\begin{array}{c}
n-3\\
\left [\frac{n-3}{2}\right ]
\end{array}
\right)+1) \cdots (\left(
\begin{array}{c}
2\\
1
\end{array}
\right)+1)\\
\leq (2^{n-3})(2^{n-4})\cdots (2)(2)\\
=2^{(n-2)(n-3)/2+1}.
\end{array}
$$

The bound $2^{(n-2)(n-3)/2+1}$ is already much lower than
doubly-exponential cost of CAD $^{[2,9]}$. We have written a
function in Maple to implement the algorithm COPOMATRIX. For
non-commercial request, we will offer for free. Please sent e-mail
to the address

\quad \quad \ yaoyong@casit.ac.cn,\ \ or,

\quad \quad  \ j.jia.xu@gmail.com.

\section{Acknowledgement}
The work of the authors were supported by the Chinese National
Science Foundation under contracts 11001228 and 10901116. The
authors also would like to thank the referees for their helpful
suggestions.





\bibliographystyle{model1-num-names}
\bibliography{<your-bib-database>}



\end{document}